\documentclass[letterpaper,11pt,reqno]{amsart}

\usepackage{graphicx,color,epsfig}
\usepackage{amsfonts,amsmath,amssymb,amsthm}

\usepackage{enumerate}

\newtheorem{thm}{Theorem}[section]
\newtheorem{prop}[thm]{Proposition}
\newtheorem{lem}[thm]{Lemma}
\newtheorem{cor}[thm]{Corollary}

\newtheorem{asm}{Assumption}

\theoremstyle{remark}
\newtheorem{rem}[thm]{Remark}

\theoremstyle{definition}

\newcommand{\ra}{\rightarrow}
\newcommand{\Ra}{\Rightarrow}

\newcommand{\R}{\mathbb R}     
\newcommand{\Z}{\mathbb Z}     

\renewcommand{\a}{\alpha}
\renewcommand{\b}{\beta}

\renewcommand{\d}{\delta}

\newcommand{\e}{\varepsilon}

\renewcommand{\l}{\lambda}

\newcommand{\s}{\sigma}

\newcommand{\iid}{i.i.d.\ }          

\newcommand{\hB}{\widehat{B}}

\newcommand{\be}{\begin{equation}}
\newcommand{\ee}{\end{equation}}

\title[Contact process with vertex-dependent infection rates]{The contact process on the complete graph with random vertex-dependent infection rates}
\author{Jonathon Peterson}
\address{Jonathon Peterson \\  Cornell University \\ Department of Mathematics \\ Malott Hall \\ Ithaca, NY 14853 \\ USA}
\email{peterson@math.cornell.edu}
\urladdr{http://www.math.cornell.edu/~peterson}
\thanks{J. Peterson was partially supported by National Science Foundation grant DMS-0802942.}

\date{May 3, 2010}


\begin{document}

\begin{abstract}
We study the contact process on the complete graph on $n$ vertices where the rate at which the infection travels along the edge connecting vertices $i$ and $j$ is equal to $\l w_i w_j / n$ for some $\l>0$, where $w_i$ are i.i.d.\ vertex weights. We show that when $E[ w_1^2] < \infty$ there is a phase transition at $\l_c > 0$ so that for $\l<\l_c$ the contact process dies out in logarithmic time, and for $\l>\l_c$ the contact process lives for an exponential amount of time. Moreover, we give a formula for $\l_c$ and when $\l>\l_c$ we are able to give precise approximations for the probability a given vertex is infected in the quasi-stationary distribution.

Our results are consistent with a non-rigorous mean-field analysis of the model. This is in contrast to some recent results for the contact process on power law random graphs where the mean-field calculations suggested that $\l_c>0$ when in fact $\l_c = 0$.
\end{abstract}

\maketitle

\section{Introduction}

The contact process is a simple model for the spread of a disease.
The standard model of the contact process on a graph $G=(V,E)$ is described informally as follows. 
Fix a parameter $\l>0$ and a set of vertices $A\subset V$. At time $t=0$ only the vertices in $A$ are infected. As time progresses, each uninfected vertex $x$ becomes infected at rate equal to $\l$ times the number of currently infected neighbors, and each infected vertex becomes a healthy (uninfected) vertex at rate 1. 
More formally, the contact process is a continuous time Markov process $\eta_t \in \{0,1\}^V$ with generator
\[
 \mathcal{L}f(\eta) = \sum_{x\in V} \left( \eta(x) + (1-\eta(x)) \l \sum_{y\sim x} \eta(y) \right) \left( f(\eta^x) - f(\eta) \right),
\]
where $f:\{0,1\}^V \ra \R$ is a bounded function, $x\sim y$ means that vertices $x$ and $y$ are connected by an edge in $E$, and $\eta^x$ is the configuration obtained from $\eta$ by switching the value of $\eta(x)$. That is, 
\be\label{switchstate}
 \eta^x(z) = \begin{cases} \eta(z) & z \neq x \\ 1- \eta(z) & z = x. \end{cases}
\ee
The contact process is also sometimes referred to as the susceptible-infected-susceptible (SIS) epidemic model. 

The behavior of the contact process depends on the parameter $\l$, and as $\l$ increases the infection spreads faster and it takes a longer amount of time for the contact process to die out (i.e., reach the absorbing state of all healthy vertices). It is then natural to ask if there is a critical values of $\l$ for which the contact process exhibits a phase transition. 
The contact process on the integer lattice $\Z^d$ has been well studied, and it is known that there is a $\l_c > 0$ such that the contact process started with a single vertex infected dies out with probability one if $\l < \l_c$ and survives forever with positive probability if $\l > \l_c$ \cite{lCVEP}. 
On any finite graph the contact process always eventually dies out, and thus it is not immediately clear how to define a phase transition. 
However, for the contact process on $[-n,n]^d \subset \Z^d$ it is known that for $\l<\l_c$ (where $\l_c$ is the critical parameter for the contact process on $\Z^d$) the contact process dies out by time $C\log n$ with high probability, whereas for $\l > \l_c$ the contact process survives for time $\exp\{cn^d\}$ with positive probability \cite{lCVEP}.
In general, one says that the contact process on a family of finite graphs is sub-critical if the time until the infection dies out is logarithmic in the number of vertices and is super-critical if with positive probability the infection survives for a time that is larger than any polynomial in the number of vertices of the graph. The critical value $\l_c$ then identifies the phase transition of the contact process from sub-critical to super-critical.


The contact process has also been studied on graphs other than $\Z^d$. Probably the first such work was done by Pemantle on the contact process on infinite trees \cite{pCPTrees}. The contact process has also been studied on certain non-homogeneous classes of graphs. 
Recently, Chatterjee and Durrett \cite{cdCPRG} and Berger, Borgs, Chayes and Saberi \cite{bbcsCPRG} considered the contact process on two different models of power-law random graphs. A power-law random graph is a general term denoting a class of graphs where the distribution of the degree of a typical vertex has tails that decay like $C k^{-\a-1}$ as $k\ra\infty$ for some $\a>1$ and $C>0$.  Physicists had previously studied the contact process on power-law random graphs using non-rigorous mean-field calculations and concluded that if $\a>3$ then there was a critical value $\l_c>0$ identifying a phase transition \cite{pvMeanField1,pvMeanField2}.
However, for the two types of power law random graphs studied in \cite{bbcsCPRG} and \cite{cdCPRG} it was shown that in fact $\l_c = 0$ (i.e., the contact process survives for a long time for any $\l>0$). The long time survival of the contact process implies the existence of a quasi-stationary distribution.
The mean-field calculations suggest that the average density $\rho(\l)$ of infected sites in the quasi-stationary distribution satisfies
 $\rho(\l) \sim C \l^\b$ for some $\b>0$ as $\l\ra \l_c^+$. However, upper and lower bounds on $\rho(\l)$ calculated in \cite{cdCPRG} show that the ``critical exponent'' $\beta$ must be different from the mean-field predictions. 

In this paper, instead of studying the contact process on a random graph, we use a deterministic graph (the complete graph $K_n$) and instead make the rates at which the infections travel along edges to be random. We will choose the random infection rates in a manner that is inspired by the power-law random graph model of Chung and Lu \cite{clADRG}.
The mean field calculations in this model are exactly the same as in the case of the contact process on power law random graphs. However, in this model the mean field calculations actually turn out to be correct,
and we are able to obtain a formula for $\l_c$ and approximations for the average density $\rho(\l)$ of infected sites.

\subsection{Description of the Model}
We now turn to a description of the specific process studied in this paper. 
We will denote the set of vertices of $K_n$ by $[n]:=\{1,2,\ldots n\}$.
Let $\mathbf{w} = \{w_i\}_{i=1}^\infty $ be a sequence of non-negative real numbers, and assign weights $\mathbf{w}_n = \{w_i\}_{i=1}^n$ to the vertices of $K_n$. 
Then, given the vertex weights, the rate at which infections are transmitted from $i$ to $j$ (or $j$ to $i$) is equal to $\l w_i w_j/n$. The rate at which infected vertices become healthy is kept constant at $1$. Formally, for fixed $\mathbf{w}_n$, $\l>0$, and $A\subset [n]$ we let $\eta_t$ be the Markov process with $\eta_0(\cdot) = \mathbf{1}_{A}(\cdot)$ and generator
\[
 (\mathcal{L}_{\mathbf{w}_n} f )(\eta) = \sum_{i=1}^n \left( \eta(i) + (1-\eta(i))\sum_{j\neq i} \eta(j) \frac{\l w_i w_j}{n} \right)( f(\eta^i) - f(\eta) ),
\]
where $\eta^i$ is defined as in \eqref{switchstate}. 
The law of $\eta_t$ in this case will be denoted by $P_{\mathbf{w}}^{n,A}$. Corresponding expectations will be denoted by $E_\mathbf{w}^{n,A}$. 
Often we will be interested in the cases where either all vertices are initially infected or just a single vertex is initially infected. Thus, for notational convenience we will abbreviate $P_{\mathbf{w}}^{n,[n]}$ and $P_{\mathbf{w}}^{n,\{i\}}$ by $P_{\mathbf{w}}^{n}$ and $P_{\mathbf{w}}^{n,i}$, respectively. Corresponding expectations will be denoted similarly. 

The Markov process $\eta_t$ takes values in $\{0,1\}^n$, and so we can identify $\eta_t$ with a vector of length $n$. 
We will use $\mathbf{0}$ and $\mathbf{1}$ to denote the vectors of all zeros and all ones so that $\eta_t = \mathbf{0}$ denotes all vertices being healthy at time $t$ and $\eta_t = \mathbf{1}$ denotes all vertices being infected at time $t$.

We will make the following assumption on the vertex weights
\begin{asm}
 The sequence of vertex weights $\mathbf{w} = \{ w_i \}_{i=1}^\infty$ is an \iid sequence of random variables with common distribution $\mu$. 
\end{asm}
Expectations with respect to the measure $\mu$ on the vertex weights will be denoted $E_\mu$. 
\begin{rem}
 The essential property that we use is that 
$\lim_{n\ra\infty} \frac{1}{n} \sum_{i=1}^n f(w_i) = E_\mu[f(w_1)]$ for certain rational functions $f(w)$. Thus, the statements of the theorems and the proofs remain essentially unchanged by letting $\mathbf{w}$ be an ergodic sequence or even a ``nice'' deterministic sequence (c.f. \cite{mrDetSeq2}) with empirical distributions approximating the measure $\mu$. 
\end{rem}

Before stating our main results, we explain briefly the motivation behind this model. 
One way of constructing the contact process $\eta_t$ with random edge weights described above is to construct for each $i$ a Poisson point process $N_i(t)$ with rate 1 and for each pair $i\neq j$ a Poisson point process $N_{i,j}(t)$ with rate $\l w_i w_j/n$, such that all the Poisson point processes are independent. At each jump time of the process $N_i(\cdot)$, if the vertex $i$ is infected it becomes healthy, and at each jump time of the process $N_{i,j}(\cdot)$ if exactly one of the vertices $i$ or $j$ is infected then the other vertex becomes infected as well. 

If we consider the vertices to represent computers or individuals in a network then the jump times of the process $N_{i,j}(\cdot)$ represent connections made between the respective individuals or computers (such as human contact or e-mail message sent). 
Thus, if we only keep track of connections formed over a short time period $[0,\d]$, the resulting random graph will have the edge between $i$ and $j$ present with probability $1-\exp\{-\d \l w_i w_j/n\} \approx \d \l w_i w_j/n$. Chung and Lu \cite{clADRG} studied a model for random graphs where, given a sequence of vertex weights $w_i$, the 
probability that there is an edge connecting vertices $i$ and $j$ is proportional to $w_i w_j/n$. If the vertex weights $w_i$ have power-law tails, 
then the degree distribution of the resulting random graph has power law tails. 
Thus, if the distribution $\mu$ on vertex weights has power law tails one might expect the contact process on $K_n$ with random infection rates to be similar to the contact process on a power law random graph. 
However, our main results demonstrate that this is not the case, and that instead the mean-field predictions are actually correct.

\subsection{Main Results}

We now turn to the statements of the main results of the paper. 
We introduce these via the mean-field heuristics which help to explain them. 
Assume for now that there are only a finite number of weights and also that the expected number of vertices of weight $x$ is exactly $\mu(\{x\}) n$. 
Let $p_t(x)$ be the probability that a vertex of weight $x$ is infected at time $t$, and let $N_t(x) = p_t(x) \mu(\{x\}) n$ be the expected number of vertices of weight $x$ that are infected at time $t$.
Then, under the mean field assumption that $\eta_t(i)$ and $\eta_t(j)$ are independent (which is not true), we obtain that
\[
 \frac{d N_t(x)}{dt} = - N_t(x) + \sum_y (\mu(\{x\})n - N_t(x)) N_t(y) \frac{\l x y}{n}.
\]
Recalling that $p_t(x) = N_t(x)/(\mu(\{x\}) n)$ we obtain that
\begin{align}
 \frac{d p_t(x)}{dt} &= - p_t(x) + \sum_y (1 - p_t(x)) p_t(y) \l x y \, \mu(\{y\}) \nonumber \\
&= - p_t(x) + \int  (1 - p_t(x)) p_t(y) \l x y \, \mu(dy). \label{meanfield}
\end{align}
The above mean-field equation should also hold when the distribution of vertex weights is continuous as well. 

If the contact process survives for a long time then there should be a quasi-stationary distribution. Thus, we look for a stationary solution to \eqref{meanfield}. That is, we want to find a function $p(x)$ such that 
\be\label{mf1}
 p(x) = \int  (1-p(x))p(y) \l x y \, \mu(dy) = \int  p(y) \l x y \, \mu(dy) - p(x) \int  p(y) \l x y  \, \mu(dy). 
\ee
Solving for $p(x)$ we obtain that
\be\label{mf2}
 p(x) = \frac{\s x}{1+\s x} , \quad \text{where } \s = \l \int  y p(y) \, \mu(dy).
\ee
Substituting the formula for $p(x)$ on the left into the equation on the right yields the equation
\begin{equation}\label{mfseq}
 1 = \l \int  \frac{ y^2}{1+\s y} \, \mu(dy).
\end{equation}
We wish to characterize for what values of $\l$ (depending on the distribution $\mu$ of vertex weights) there is a $\s > 0$ which solves \eqref{mfseq}. 
Note that the right hand side of \eqref{mfseq} is decreasing in $\s$. Thus if there is a solution to \eqref{mfseq} it is unique. Moreover, (assuming $E_\mu w_1 < \infty$)
\[
 \lim_{\s \ra \infty}  \l \int  \frac{y^2}{1+\s y} \, \mu(dy) = 0
\]
and
\[
 \lim_{\s \ra 0^+} \int  \l \frac{y^2}{1+\s y} \, \mu(dy) = \int  \l y^2 \, \mu(dy) = \l E_\mu( w_1^2 ) . 
\]
This leads us to the following definitions. 
Let
\be\label{lcdef}
 \l_c := \begin{cases} \frac{1}{E_\mu w_1^2} & \text{if } E_\mu w_1^2 < \infty\\ 0 & \text{otherwise,} \end{cases}
\ee
and for $\l>\l_c$ define $\s(\l)$ as follows.  
\be\label{sldef}
 \s(\l) \text{ is the unique } \s > 0 \text{ that solves } \qquad 1 = \l E_\mu\left[ \frac{ w_1^2}{1+ \s w_1} \right]. 
\ee

Our main results are a confirmation of the above mean-field heuristics. The first result verifies the existence of a phase transition at $\l_c$. 
\begin{thm}\label{mainthm}
Let the vertex weights $\mathbf{w}$ be \iid with distribution $\mu$, and let $\l_c$ be defined as in \eqref{lcdef}. 
\begin{enumerate}[(i)]
 \item \label{dieout} If $\l < \l_c$, then there exists a constant $C>0$ such that 
\[
 \lim_{n\ra\infty} P_\mathbf{w}^n( \eta_{C\log n} = \mathbf{0} ) = 1, \quad \mu-a.s.
\]
 \item \label{mainthmepidemic} If $\l> \l_c$, then there exists a constant $c>0$ such that 
\[
 \lim_{n\ra\infty} P_\mathbf{w}^n( \eta_{e^{cn}} \neq \mathbf{0} ) = 1, \quad \mu-a.s.
\]
\end{enumerate}
\end{thm}
\begin{rem}
 When $E_\mu w_1^2 = \infty$, then $\l_c=0$ and condition \eqref{dieout} is vacuous. That is, the contact process is always super-critical if $E_\mu w_1^2 = \infty$. 
\end{rem}
\begin{rem}
 The probabilities in the conclusion of Theorem \ref{mainthm} are actually random variables, since they depend on the vertex weights $\mathbf{w}$ which are random. 
 The conclusion of Theorem \ref{mainthm} is that the above limits hold for $\mu$-almost every realization of the vertex weights $\mathbf{w}$.
Such limits are often called \emph{quenched} limiting statments in the field of random environments. 
\end{rem}

Our second main result states that when the contact process is super-critical, the probability of a vertex being infected may be uniformly approximated by the mean-field prediction \eqref{mf2}. 
\begin{thm}\label{mainthm2}
 Let the vertex weights $\mathbf{w}$ be \iid with distribution $\mu$, and let $\l_c$ and $\s(\l)$ be defined as in \eqref{lcdef} and \eqref{sldef}, respectively. 
Then, for any $\l>\l_c$ and $\e> 0$ there exist constants $C,c>0$ (depending on $\e$) so that 
\[
 \limsup_{n\ra\infty} \sup_{t\in[C\log n, e^{cn}]} \sup_{i\in [n]} \left| P_\mathbf{w}^n( \eta_t(i) = 1 ) - \frac{\s(\l) w_i}{1+\s(\l)w_i}  \right| \leq \e, \quad \mu-a.s. 
\]
\end{thm}
\begin{rem}
 We are actually able to prove an upper bound on the probability of a vertex being infected that is slightly better than what is implied by the statement of Theorem \ref{mainthm2}. See Proposition \ref{rhoub} for a precise statement. 
\end{rem}

In their study of the contact process on a random graph, Chatterjee and Durrett \cite{cdCPRG} analyzed the contact process at time $t=e^{\sqrt{n}}$. 
The time $e^{\sqrt{n}}$ is large enough for the contact process to have stabilized but small enough so that with high probability it has not died out. Thus, the distribution of the contact process at time $t=e^{\sqrt{n}}$ when started with all vertices infected was called the quasi-stationary distribution. 
As in \cite{cdCPRG}, we define $\rho_n(\l)$ to be the expected number of infected vertices at time $e^{\sqrt{n}}$ when initially all vertices are infected. That is,
\[
 \rho_n(\l) = E_\mu\left[ \frac{1}{n} \sum_{i=1}^n P_\mathbf{w}^{n}(\eta_{e^{\sqrt{n}}}(i) = 1 ) \right].
\]
Then, we obtain the following simple Corollary of Theorem \ref{mainthm2}. 
\begin{cor}\label{rholim}
If $\l > \l_c$, then 
 \[
  \rho(\l) := \lim_{n\ra\infty} \rho_n(\l) = E_\mu \left[ \frac{\s(\l) w_1}{1+\s(\l)w_1}  \right].
 \]
\end{cor}
As mentioned above, upper and lower bounds on the critical exponent of $\rho_n(\l)$ were derived in \cite{cdCPRG}. 
Since
\[
 E_\mu \left[ \frac{\s w_1}{1+\s w_1}  \right] \sim E_\mu[ w_1] \s, \quad \text{ as } \s \ra 0^+,
\]
Corollary \ref{rholim} implies that $\rho(\l)$ and $\s(\l)$ have the same critical exponent. 
The following Proposition allows us to not only identify the critical exponent of $\s(\l)$ but also the leading constants. 

\begin{prop}\label{cexp}
 Let $\mu(w_1 > x) \sim C x^{-(\a-1)}$ as $x\ra\infty$ for some $C>0$ and $\a>2$, and let $\l_c$ and $\s(\l)$ be defined as in \eqref{lcdef} and \eqref{sldef}. Then, as $\d\ra 0^+$,
\begin{align*}
\s(\d) &\sim \left( \frac{C \pi (\a-1)}{\sin(\pi \a)} \right)^{1/(3-\a)} \d^{1/(3-\a)} & \text{ if } \a\in(2,3)\\
\log \s(\d) &\sim \frac{-1}{2C\d} & \text{ if } \a = 3\\
\s(\l_c + \d) &\sim \left( \frac{-\sin(\pi \a)}{C \l_c^2 (\a-1)\pi} \right)^{1/(\a-3)} \d^{1/(\a-3)} & \text{ if } \a\in(3,4)\\
\s(\l_c + \d) &\sim  \frac{1}{3 C \l_c^2} \frac{\d}{\log(1/\d)} & \text{ if } \a=4\\
\s(\l_c + \d) &\sim \frac{1}{\l_c^2 E_\mu[w_1^3]} \d & \text{ if } \a>4.
\end{align*}
Moreover, the conclusion in the case $\a>4$ also holds under the assumption that $E_\mu[w_1^3] < \infty$ (i.e., without any assumption on tail asymptotics of $\mu$).  
\end{prop}
\begin{rem}
 Since Theorem \ref{mainthm2} shows that the infection probabilities agree with the mean-field predictions, these critical exponents agree with the mean-field predictions in \cite{pvMeanField2}. However, the case $\a=4$ was not considered in \cite{pvMeanField2}. 
\end{rem}
\begin{rem}
 Proposition \ref{cexp} states that in the case $\a=3$, $\s(\d) = \exp\{ -\frac{1}{\d}(\frac{1}{2C} + o(1)) \}$. Riordan \cite{rSGCSFRG} has previously observed a similar fast rate of decay for the fraction of vertices in the giant component of percolation on the Barb\'asi-Albert model of power law random graphs (which have $\a=3$). 
\end{rem}

Before proceeding to the proofs of the above results, we recall two useful facts about the contact process that we will use. Both of these facts are found in \cite{lCVEP} and can be proved using what is called the ``graphical representation.''
\begin{itemize}
 \item \textbf{Monotonicity in infection rates.} Increasing the infection rate along any edge only increases the number of vertices infected at any time $t$. In particular, if $\mathbf{\widetilde{w}} \leq \mathbf{w}$ in the sense that $\widetilde{w}_i \leq w_i$ for all $i$, an easy coupling argument implies that
\be\label{monotonicity}
 \mathbf{\widetilde{w}} \leq \mathbf{w} \Ra P_\mathbf{\widetilde{w}}^{n,A}( \eta_t \geq \mathbf{1}_B ) \leq P_\mathbf{w}^{n,A}( \eta_t \geq \mathbf{1}_B ), \quad \forall A,B\subset[n] .  
\ee
 \item \textbf{Self-duality.} For any subsets of vertices $A$ and $B$, the probability that an element of $B$ is infected at time $t$ when the infection starts from $A$ is equal to the probability that an element of $B$ is infected at time $t$ when the infection starts from $B$. That is, 
\[
 P_{\mathbf{w}}^{n,A}( \eta_t(i) = 1, \, \text{ for some } i\in B ) = P_{\mathbf{w}}^{n,B}( \eta_t(i) = 1, \, \text{ for some } i\in A ).
\]
An important special case of this is when $A=\{i\}$ and $B = [n]$. In this case, we obtain that
\be\label{duality}
 P_\mathbf{w}^{n,i}(\eta_t \neq \mathbf{0}) = P_\mathbf{w}^n( \eta_t(i) = 1). 
\ee
\end{itemize}

The remainder of the paper is organized as follows. In Section \ref{mdbirthdeath} we consider the special case when the distribution $\mu$ has finite support. In this case we identify the contact process with a multi-dimensional birth-death chain. By analyzing the birth-death chain we then show that when $\l>\l_c$ the contact process survives for time $e^{cn}$ with high probability, and we also obtain a lower bound on the probability a given vertex is infected at any time $t\leq e^{cn}$. In Section \ref{finitetocontsupport} we extend the results of Section \ref{mdbirthdeath} to the general case by approximating $\mu$ by a measure with finite support. 
Next in Section \ref{multitypeBP} we complete the proofs of Theorems \ref{mainthm} and \ref{mainthm2} by proving the complementary results to those obtained in Sections \ref{mdbirthdeath} and \ref{finitetocontsupport}. That is, we show that the contact process dies out with high probability by time $C \log n$ when $\l < \l_c$, and we obtain upper bounds on the probability a given vertex is infected at any time $t \geq C \log n$ when $\l > \l_c$. The main technique used in Section \ref{multitypeBP} is a comparison of the contact process with a related multi-type branching process. 
Finally, the proof of Proposition \ref{cexp} is given in Section \ref{criticalexponent}.

\textbf{Acknowledgement.} Many thanks to Rick Durrett for the many helpful discussions in the process of writing this paper.

\section{Finitely many types}\label{mdbirthdeath}
In this section we will analyze the above model under the assumption that the support of the $\mu$ on vertex weights $w_i$ is finite. Let
$\mathbf{W} = (W_1, W_2, \ldots, W_m)$ be the possible vertex weights, each occuring with probability $p_i = \mu(w_1 = W_i)$. 
We will classify the vertices according to the weight that they are assigned. That is, a vertex with weight $W_i$ will be referred to as a \emph{type $W_i$} vertex. 
Since all type $W_i$ vertices are equivalent, we need only to keep track of the number $X_i(t)$ of type $W_i$ vertices that are infected at time $t$. 
That is,
\[
 X_i(t):= \# \left\{ j : w_j = W_i, \text{ and } \eta_t(j) = 1 \right\}. 
\]
Note that $\mathbf{X}(t) = ( X_1(t), X_2(t), \ldots, X_m(t))$ is an $m$-dimensional birth death process.
Given the vertex weights $\mathbf{w}$, we denote the law of $\mathbf{X}(\cdot)$ started at $\mathbf{X}(0) = \mathbf{x}$ by
$P_\mathbf{w}^{n,\mathbf{x}}$. As with the contact process, we will write $P_\mathbf{w}^n$ for the law of $\mathbf{X}(\cdot)$ when the associated contact process starts with all vertices infected (that is, when $\sum_{j=1}^m X_j(0) = n$).

The main result of this section is the following, which not only shows that the contact process is super-critical when $\l>\l_c$, but also gives a lower bound on the probability of a given vertex to be infected. 

\begin{prop}\label{finitesupport}
 Let the distribution $\mu$ have finite support, and let $\l_c$ and $\s(\l)$ be defined as in \eqref{lcdef} and \eqref{sldef}. Then, if $\l > \l_c$ and  $\eta \in (0, \s(\l))$ there exists a constant $C>0$ such that
\[
 \lim_{n\ra\infty} P_\mathbf{w}^n \left( \inf_{t\leq e^{Cn}} X_i(t) \geq \frac{\eta W_i}{1+\eta W_i} p_i n, \, \forall i=1,2,\ldots, m \right) = 1, \quad \mu-a.s.
\]
\end{prop}

\begin{proof}
 For any $\e>0$, let 
\[
 U_n(\e) := \left \{ \mathbf{x}\in\Z^m : x_i \geq \frac{\e W_i}{1+\e W_i} p_i n\, , \, \forall i=1,2,\ldots, m \right\}. 
\]
Also, let $\mathbf{X}[t_1,t_2] := \{ \mathbf{X}(t): t \in [t_1,t_2] \}$ be the trace of the process $\mathbf{X}(t)$ between times $t_1$ and $t_2$. 
The idea of the proof is that there exist constants $\d \in(\eta,\s(\l))$, $\tau>0$, and $C'>0$ that are independent of $n$ and such that $\mu-a.s.$,
\be\label{staydontexit}
 \inf_{\mathbf{x} \in U_n(\d)} P_\mathbf{w}^{n,\mathbf{x}}\left( \mathbf{X}(\tau) \in U_n(\d) \, , \,  \mathbf{X}[0,\tau]\subset U_n(\eta) \right) \geq 1 - e^{-C' n}, \quad \text{for all } n \text{ large enough.}
\ee
That is, we can (for $n$ sufficiently large) uniformly bound from below the probability that starting from a point $\mathbf{x}\in U_n(\d)$ the process $\mathbf{X}(t)$ a short time later is still in $U_n(\d)$ and hasn't exited $U_n(\eta)$. 
By dividing $[0,e^{C'n/2}]$ into $e^{C'n/2}/\tau$ intervals of length $\tau$, \eqref{staydontexit} implies that, $\mu-a.s.$, 
\[
\lim_{n\ra\infty} P_\mathbf{w}^n \left( \mathbf{X}[0,e^{C'n/2}] \subset U_n(\eta)  \right) \geq \lim_{n\ra\infty} 1-e^{-C'n/2}/\tau = 1.
\]
This last statement is equivalent to the conclusion of the proposition. 

It remains only to prove \eqref{staydontexit}. 
The monotonicity of the contact process implies that we only need a lower bound on the inner probability in \eqref{staydontexit} at a single point. That is, 
\begin{align*}
\inf_{\mathbf{x} \in U_n(\d)} P_\mathbf{w}^{n,\mathbf{x}}\left( \mathbf{X}(\tau) \in U_n(\d) \, , \,  \mathbf{X}[0,\tau] \subset U_n(\eta) \right) = P_\mathbf{w}^{n,\mathbf{x}^\d}\left( \mathbf{X}(\tau) \in U_n(\d)\, , \,  \mathbf{X}[0,\tau] \subset U_n(\eta) \right),
\end{align*}
where
\[
 \mathbf{x}^{\d} = (x^\d_1,x^\d_2,\ldots,x^\d_m),\qquad x^\d_i := \left\lceil  \frac{\d W_i}{1+\d W_i} p_i n \right\rceil. 
\]
For any $0<\d<\e$ let 
\[
 B_n(\d,\e):= \left\{ \mathbf{x} \in \Z^m: \frac{\d W_i}{1+\d W_i} \leq \frac{x_i}{p_i n} \leq \frac{\e W_i}{1+\e W_i}, \; i=1,2,\ldots m \right\}.
\]
Note that $B_n(\d,\e) \subset U_n(\d)$. 
Then, to prove \eqref{staydontexit} it is enough to show that for some $\eta<\d<\e<\s(\l)$, and constants $\tau,C'>0$, $\mu-a.s.$,
\be\label{staydontexit2}
 P_\mathbf{w}^{n,\mathbf{x}^\d}\left( \mathbf{X}(\tau) \in U_n(\d) \, , \,  \mathbf{X}[0,\tau] \subset B_n(\eta,\e) \right) \geq 1 - e^{-C' n},   \quad \text{for all } n \text{ large enough.}
\ee
The idea behind proving \eqref{staydontexit2} is that by choosing $\d$ close enough to $\eta < \s(\l)$, the drift of the process $\mathbf{X}(t)$ is increasing in all coordinates in a neighborhood of $\mathbf{x}^\d$ (with high probability). Then, the process $\mathbf{X}(t)$ starting at $\mathbf{x}^{\d}$ will with high probability increase in all coordinates by time $\tau$ and by choosing $\tau>0$ small enough it will also not have exited $B_n(\eta,\e)$.

\begin{figure}
 \includegraphics{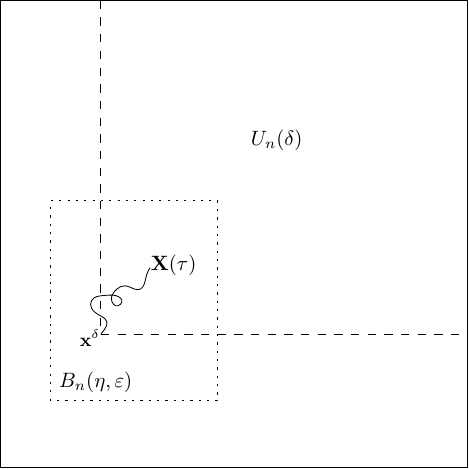}
\caption{An example of the event in \eqref{staydontexit2}. The key is to choose $\eta<\d<\e$ so that the drift of the process $\mathbf{X}(t)$ is increasing in all coordinates inside of $B_n(\eta,\e)$.}
\end{figure}

As mentioned above, $\mathbf{X}(t)$ is an $m$-dimensional birth death process (given the vertex weights), but in order to describe the jump rates we need some additional notation. 
Let 
\[
 N_{i,n} := \# \{ j\leq n: w_j = W_i \}, \quad \text{and} \quad \hat{p}_{i,n} := \frac{N_{i,n}}{n}
\]
be the number and proportion, respectively, of type $W_i$ vertices. Note that the law of large numbers implies that
\[
 \lim_{n\ra\infty} \hat{p}_{i,n} = p_i, \quad \mu-a.s.
\]
Let $\{\mathbf{e}_i\}_{i=1}^m$ be the canonical basis for $\Z^m$. If $\mathbf{x}=(x_1,x_2,\ldots,x_m) \in \Z^m_+$, then when $\mathbf{X}(t)=\mathbf{x}$ the process increases in the $i$th coordinate to $\mathbf{x} + \mathbf{e}_i$ at rate $q_i^+(\mathbf{x})$ and decreases in the $i$th coordinate to $\mathbf{x} - \mathbf{e}_i$ at rate $q_i^-(\mathbf{x})$, where 
\be\label{jumpratesq}
 q_i^+(\mathbf{x}) = (N_{i,n} - x_i)\sum_{j=1}^m x_j \frac{\l W_i W_j}{n}, \quad \text{and}\quad  q_i^-(\mathbf{x}) = x_i. 
\ee
The following Lemma allows us to bound the jump rates in some $B_n(\eta,\e)$. 
\begin{lem}\label{boundedrates}
For any $\eta < \s(\l)$, there exists an $\e \in (\eta,\s(\l))$ and positive constants $\{\a_i^+,\a_i^-\}_{i=1}^m$ such that $\mu-a.s.$, the event
\[
 \bigcap_{i=1}^m \bigcap_{\mathbf{x} \in B_n(\eta,\e)} \left\{ q^-_i(\mathbf{x}) \leq \a_i^- n < \a_i^+ n \leq  q^+_i(\mathbf{x}) \right\}
\]
occurs for all $n$ large enough. 
\end{lem}
Postponing for now the proof of Lemma \ref{boundedrates} we continue with the proof of \eqref{staydontexit2}. 
Let $\widetilde{\mathbf{X}}(t) = (\widetilde{X}_1(t),\ldots, \widetilde{X}_m(t))$, where the $\widetilde{X}_i(t)$ are independent continuous time simple random walks that increase by 1 at rate $\a^+_i n$ and decrease by 1 at rate $\a^-_i n$. Then \eqref{boundedrates} implies that (for $n$ large enough) $\mathbf{X}(t)$ stochastically dominates $\widetilde{\mathbf{X}}(t)$ while inside of $B_n(\eta,\e)$. Let $\widetilde{P}^\mathbf{x}$ be the law of $\widetilde{\mathbf{X}}(t)$ started at $\widetilde{\mathbf{X}}(0)=\mathbf{x}$. 
Since $\widetilde{\mathbf{X}}(t)$ is a continuous time simple random walk with positive drift in all coordinates and jump rates proportional to $n$, standard large deviation estimates yield the following Lemma whose proof we omit. 
\begin{lem}\label{Xtilde}
 Let $\a_i^-<\a_i^+$ for $i=1,2,\ldots, m$ be given. Then, there exists a constant $C''>0$ (depending only on the $\a_i^{\pm}$) such that for any $\d>0$,
\[
 \widetilde{P}^{\mathbf{x}^\d}\left( \widetilde{\mathbf{X}}(t) \notin U_n(\d) \right) \leq m e^{-C'' t n} , \quad \forall t\geq 0. 
\]
\end{lem}
On the event $\bigcap_{i=1}^m \bigcap_{\mathbf{x} \in B_n(\eta,\e)} \left\{  q^-_i(\mathbf{x}) \leq \a_i^- n < \a_i^+ n \leq  q^+_i(\mathbf{x}) \right\}$, since $\mathbf{X}(t)$ stochastically dominates $\widetilde{\mathbf{X}}(t)$ inside of $B_n(\eta,\e)$ we have
\be \label{XtildeX}
 \widetilde{P}^{\mathbf{x}^\d}\left( \widetilde{\mathbf{X}}(\tau) \in U_n(\d) \right) \leq P_\mathbf{w}^{n,\mathbf{x}^\d}\left( \mathbf{X}(\tau) \in U_n(\d) \, , \, \mathbf{X}[0,\tau] \subset B_n(\eta,\e) \right) 
+ P_\mathbf{w}^{n,\mathbf{x}^\d}\left( \mathbf{X}[0,\tau] \nsubseteq B_n(\eta,\e) \right)
\ee
Now, the total jump rate of the process $\mathbf{X}(t)$ an any point $\mathbf{x}$ is always bounded by $C n$ for some constant $C$ (depending on $\l$ and the $W_j$). 
Then, since the distance from $\mathbf{x}^\d$ to the complement of $B_n(\eta,\e)$ is at least $c n$ for some constant $c$ (depending on $\eta,\d,\e$ and the $W_j$) by choosing $\tau>0$ sufficiently small we obtain that for some $C_0>0$,
\[
 P_\mathbf{w}^{n,\mathbf{x}^\d}\left( \mathbf{X}[0,\tau] \nsubseteq B_n(\eta,\e) \right) \leq e^{-C_0 n}. 
\]
Thus, applying Lemmas \ref{boundedrates} and \ref{Xtilde} to \eqref{XtildeX} we obtain that $\mu-a.s.$,
\[
 P_\mathbf{w}^{n,\mathbf{x}^\d}\left( \mathbf{X}(\tau) \in U_n(\d) \, , \, \mathbf{X}[0,\tau] \subset B_n(\eta,\e) \right) \geq 1-m e^{-C'' \tau n} - e^{-C_0 n} \geq 1-e^{-C' n}, 
\]
for all large $n$ if $C' < \min\{ C'', C_0 \}$. This completes the proof of \eqref{staydontexit2} and thus also (modulo the proof of Lemma \ref{boundedrates}) the proof of Proposition \ref{finitesupport}. 
\end{proof}

We now return to the proof of Lemma \ref{boundedrates}.
\begin{proof}[Proof of Lemma \ref{boundedrates}]
We ultimately wish to control $q_i^+(\mathbf{x})$ and $q_i^-(\mathbf{x})$ inside of some $B_n(\eta,\e)$, but fluctuations in the number of type $W_i$ vertices force us to deal first with a slightly different set that depends on the actual number of type $W_i$ vertices. 
For any $0<\eta<\e$ define
\[
 \hB_n(\eta,\e) := \left\{ \mathbf{x} \in \Z^m: \frac{\eta W_i}{1+\eta W_i} \leq \frac{x_i}{N_{i,n}} \leq \frac{\e W_i}{1+\e W_i}, \; i=1,2,\ldots m \right\}.
\]
Note that the region $\hB_n(\eta,\e)$ depends on what the actual vertex weights are, whereas the region $B_n(\eta,\e)$ only depends on the given parameters. However, since $N_{i,n}/n \ra p_i$, $\mu-a.s.$ we obtain that for any $\eta'< \eta < \e< \e'$
\be\label{Lsubset}
 B_n(\eta,\e) \subset \hB_n(\eta',\e')  ,   \quad \text{for all } n \text{ large enough}, \, \mu-a.s.
\ee

Recall the formulas for the jump rates $q_i^+(\mathbf{x})$ and $q_i^-(\mathbf{x})$. Then, for any $\mathbf{x} \in \hB_n(\eta,\e)$, 
\begin{align*}
 q_i^+(\mathbf{x}) &\geq (N_{i,n} - \frac{\e W_i}{1+\e W_i} N_{i,n})\sum_{j=1}^m \frac{\eta W_j}{1+\eta W_j}N_{j,n} \frac{\l W_i W_j}{n} \\
&= \l N_{i,n}\frac{W_i}{1+\e W_i} \sum_{j=1}^m \frac{\eta W_j^2}{1+\eta W_j} \frac{N_{j,n}}{n} \\
&= \l \hat{p}_{i,n} \frac{W_i}{1+\e W_i} \sum_{j=1}^m \hat{p}_{j,n} \frac{\eta W_j^2}{1+\eta W_j}   n \\
&=: \hat{\theta}_i^+(\eta,\e,n) n 
\end{align*}
and 
\[
 q_i^-(\mathbf{x}) \leq \frac{\e W_i}{1+\e W_i} N_{i,n} = \frac{\e W_i}{1+\e W_i} \hat{p}_{i,n} n =: \hat{\theta}_i^-(\eta,\e,n) n
\]
Again, since $\hat{p}_{i,n} \ra p_i$ we obtain
\be\label{thetalim1}
 \lim_{n\ra\infty} \hat{\theta}_i^+(\eta,\e,n) = \l p_i \frac{W_i}{1+\e W_i} \sum_{j=1}^m p_j \frac{\eta W_j^2}{1+\eta W_j} =: \theta_i^+(\eta,\e), \quad \mu-a.s.
\ee
and
\be\label{thetalim2}
 \lim_{n\ra\infty} \hat{\theta}_i^-(\eta,\e,n) = p_i \frac{\e W_i}{1+\e W_i} =: \theta_i^-(\eta,\e), \quad \mu-a.s.
\ee
Note that $\theta_i^+ = (1+\Delta) \theta_i^-$ where 
\[
 \Delta = \Delta(\eta,\e) = \frac{\l \eta}{\e} \sum_{j=1}^m \frac{p_j W_j^2}{1+\eta W_j} - 1 = \frac{\l\eta}{\e} E_\mu\left[ \frac{w_1^2}{1+\eta w_1} \right] - 1. 
\]

Now, the definition of $\s(\l)$ implies that $\Delta(\eta,\eta)= \l E_\mu[ w_1^2/(1+\eta w_1) ] - 1 > 0$ for any $\eta < \s(\l)$. Therefore, we can choose an $\eta' < \eta < \e' < \s(\l)$ such that $\Delta(\eta',\e') > 0$. 
Fix such a $\eta'$ and $\e'$ and choose $\a_i^+$ and $\a_i^-$ so that
\[
 \theta_i^-(\eta',\e') < \a_i^- < \a_i^+ < \theta_i^+(\eta',\e'). 
\]
Then, \eqref{thetalim1} and \eqref{thetalim2} imply that 
\be\label{thetahatbounds}
  \hat{\theta}_i^-(\eta',\e',n) \leq \a_i^- < \a_i^+ \leq \hat{\theta}_i^+(\eta',\e',n) ,   \quad \text{for all } n \text{ large enough}, \, \mu-a.s.
\ee
Now choose $\e \in (\eta,\e')$. 
Applying \eqref{Lsubset} and \eqref{thetahatbounds} completes the proof of Lemma \ref{boundedrates}
\end{proof}

We end this section with the following Corollary. 

\begin{cor}\label{iplb}
Let the distribution $\mu$ have finite support. Then, for any $\l > \l_c$ and $\eta< \s(\l)$, there exists a constant $c>0$ such that
 \[
  \liminf_{n\ra\infty} \inf_{t\leq e^{cn}} \inf_{i\in[n]} \left( P_\mathbf{w}^n( \eta_{t}(i) = 1 ) - \frac{\eta w_i}{1+\eta w_i} \right) \geq 0,  \quad \mu-a.s.
 \]
\end{cor}
\begin{proof}
For $\l>\l_c$ and $\eta<\s(\l)$ fixed, chose $c>0$ as in the statement of Proposition \ref{finitesupport}. 
Suppose that vertex $i$ is a type $W_j$ vertex. That is $w_i = W_j$. Then, since there are $N_{j,n}$ vertices of type $W_j$, 
\[
 E_{\mathbf{w}}^n [ X_j(t) ] = N_{j,n} P_\mathbf{w}^n( \eta_{t}(i) = 1 ).  
\]
Therefore,
\begin{align}
 \inf_{t\leq e^{cn}} \inf_{i\in[n]} \left( P_\mathbf{w}^n( \eta_{t}(i) = 1 ) - \frac{\eta w_i}{1+\eta w_i} \right)
&= \inf_{t\leq e^{cn}} \inf_{j\leq m} \left( \frac{1}{N_{j,n}} E_{\mathbf{w}}^n [ X_j(t) ] - \frac{\eta W_j}{1+\eta W_j} \right) \nonumber \\
& \geq \inf_{j\leq m} \left( \frac{1}{N_{j,n}} E_{\mathbf{w}}^n \left[ \inf_{t\leq e^{cn}} X_j(t) \right] - \frac{\eta W_j}{1+\eta W_j} \right). \label{infinfX}
\end{align}
Now, Proposition \ref{finitesupport} and the fact that $N_{j,n} \sim p_{j} n$ as $n\ra\infty$ imply that
\be\label{EXlim}
 \liminf_{n\ra\infty} \frac{1}{N_{j,n}} E_{\mathbf{w}}^n \left[ \inf_{t\leq e^{cn}} X_j(t) \right] \geq \frac{\eta W_j}{1+\eta W_j}, \quad \mu-a.s.
\ee
Combining \eqref{infinfX} and \eqref{EXlim} completes the proof. 
\end{proof}

\section{Infinitely many types} \label{finitetocontsupport}
In this section we will extend the results of the previous section to the case where $\mu$ does not have finite support. 
The first main result is a proof of part \eqref{mainthmepidemic} of Theorem \ref{mainthm}. 
\begin{proof}[Proof of Theorem \ref{mainthm}, part \eqref{mainthmepidemic}]
For any integer $m\geq 1$ let 
\[
 \kappa_m(x) := \frac{\lfloor x m \rfloor}{m} \wedge m. 
\]
Note that $\kappa_m(x) \leq x$ and $\lim_{m\ra\infty} \kappa_m(x) = x$ for all $x\in\R$. 
Given the sequence of vertex weights $\mathbf{w} = \{w_i\}_{i=1}^\infty$, let
\[
 \mathbf{w}^{(m)}  := \left\{ \kappa_m(w_i) \right\}_{i=1}^\infty. 
\]
The sequence vertex weights $\mathbf{w}^{(m)}$ takes on only finitely many values, and thus we may apply the results from the previous section. 
For any vertex distribution $\mu$ and $m\geq 1$, let $\l_c^{(m)}$ and $\s^{(m)}(\l)$ be the analogues of $\l_c$ and $\s(\l)$ for the modified vertex weights $\mathbf{w}^{(m)}$. That is, 
\[
 \l_c^{(m)} := \left( E_\mu[ \kappa_m(w_1)^2 ] \right)^{-1}, 
\]
and for $\l > \l_c^{(m)}$
\[
 \s^{(m)}(\l) \text{ is the unique } \s>0 \text{ that solves } \qquad 1 = \l E_\mu \left[ \frac{\kappa_m(w_1)^2}{1+ \s \kappa_m(w_1)} \right]. 
\]
The monotone convergence theorem implies that $\l_c^{(m)} \searrow \l_c$ and $\s^{(m)}(\l) \nearrow \s(\l)$ as $m\ra\infty$.

If $\l>\l_c$, then there exists some $m$ large enough so that $\l > \l_c^{(m)}$. For this $m$ fixed, we may apply Proposition \ref{finitesupport} to the modified vertex weights $\mathbf{w}^{(m)}$ to obtain that there exists a $C>0$ such that
\be\label{msurvival}
 \lim_{n\ra\infty} P_{\mathbf{w}^{(m)}}^n \left(\eta_{e^{C n}} \neq \mathbf{0} \right) = 1, \quad \mu-a.s.
\ee
Since $\mathbf{w}^{(m)} \leq \mathbf{w}$, the monotonicity of the contact process \eqref{monotonicity} implies that $P_\mathbf{w}^n(\eta_t \neq \mathbf{0} ) \geq P_{\mathbf{w}^{(m)}}^n(\eta_t \neq \mathbf{0} )$, and thus \eqref{msurvival} holds with $\mathbf{w}$ in place of $\mathbf{w}^{(m)}$. 
\end{proof}

We can also extend Corollary \ref{iplb} to the general case. This proves the lower bound needed for the proof of Theorem \ref{mainthm2}. 

\begin{cor}\label{iplb2}
Let $\l > \l_c$. Then, for any $\e>0$ there exists a constant $c>0$ so that 
 \[
  \liminf_{n\ra\infty} \inf_{t\leq e^{cn}} \inf_{i\in[n]} \left( P_\mathbf{w}^n( \eta_{t}(i) = 1 ) - \frac{\s(\l) w_i}{1+\s(\l) w_i} \right)  \geq - \e, \quad \mu-a.s.
 \]
\end{cor}
\begin{proof}
First note that an exercise in calculus shows that 
\be\label{calculus}
 \sup_{x\geq 0} \left( \frac{\s x}{1+\s x} - \frac{\eta x}{1+\eta x} \right) \leq \frac{\s - \eta}{(\sqrt{\s} + \sqrt{\eta})^2}, \qquad \forall 0\leq \eta \leq \s. 
\ee
Also, note that $\frac{\eta x}{1+\eta x}$ has derivative bounded above by $\eta$ on $x\geq 0$, and $x-\kappa_m(x)\leq 1/m$. Thus, by considering separately the cases $x\leq m$ and $x\geq m$
\[
 \sup_{x\geq 0} \left( \frac{\eta x}{1+\eta x} - \frac{\eta \kappa_m(x)}{1+\eta \kappa_m(x)} \right) \leq \max \left\{ \frac{\eta}{m} , \frac{1}{1+\eta m} \right\}, \quad \forall \eta\geq 0, \; \forall m \geq 1.
\]
Therefore, given $\l>\l_c$ and $\e>0$ we may choose $\eta<\s(\l)$ such that for all integers $m$ large enough
\be\label{etam}
 \sup_{x\geq 0} \left( \frac{\s(\l) x}{1+\s(\l) x} - \frac{\eta \kappa_m(x)}{1+\eta \kappa_m(x)} \right) < \e. 
\ee

For this choice of $\eta$, choose $m$ large enough so that $\l > \l_c^{(m)}$, $\eta < \s^{(m)}(\l)$, and \eqref{etam} holds. Then, there exists a constant $c>0$ such that the conclusion of Corollary \ref{iplb} holds for this $\eta$ and the modified vertex sequence $\mathbf{w}^{(m)}$. 
Then, the monotonicity property \eqref{monotonicity}, the choice of $\eta$ and $m$ satisfying \eqref{etam},  and Corollary \ref{iplb} imply 
\begin{align*}
& \liminf_{n\ra\infty} \inf_{t\leq e^{cn}} \inf_{i\in[n]} \left( P_\mathbf{w}^n( \eta_{t}(i) = 1 ) - \frac{\s(\l) w_i}{1+\s(\l) w_i} \right)  \\
&\qquad \geq \liminf_{n\ra\infty} \inf_{t\leq e^{cn}} \inf_{i\in[n]} \left( P_{\mathbf{w}^{(m)}}^n( \eta_{t}(i) = 1 ) - \frac{\eta \kappa_m(w_i)}{1+\eta \kappa_m(w_i)} - \e \right)  \geq -\e.
\end{align*}
\end{proof}

\section{Comparison with a multi-type branching process} \label{multitypeBP}
In this section we introduce a related multi-type branching process (MTBP) that we will use to show that the contact process dies out quickly if $\l < \l_c$. 
For fixed vertex weights $\mathbf{w}$ and $n\geq 1$, let $\mathbf{Z}(t) = (Z_1(t), Z_2(t), \ldots, Z_n(t))$ be a MTBP with $n$ types. Each individual dies at rate $1$, and individuals of type $i$ give birth to individuals of type $j$ at rate $\l w_i w_j/n$. 
The relevance of the MTBP $\mathbf{Z}(t)$ to the contact process $\eta_t$ is given by the following lemma. 
\begin{lem}\label{MTBPcouple}
 The MTBP $\mathbf{Z}(t)$ stochastically dominates the contact process $\eta_t$. That is, $\mathbf{Z}(t)$ and $\eta_t$ can be coupled in such a way that $\eta_0 = \mathbf{Z}(0)$ and $\eta_t \leq \mathbf{Z}(t)$ for all $t\geq 0$. 
\end{lem}
\begin{proof}
 We describe the coupling informally as follows. We will divide the population of the branching process into two groups, one of which will be identified with the contact process $\eta_t$. At any time $t$, we will denote the individuals in groups 1 and 2 by $\mathbf{Z}^{(1)}(t)$ and $\mathbf{Z}^{(2)}(t)$, respectively, so that $\mathbf{Z}(0) = \mathbf{Z}^{(1)}(0) + \mathbf{Z}^{(2)}(0)$. 
Let $\mathbf{Z}^{(1)}(0) = \eta_0$, so that initially all individuals are in group 1. The dynamics of the processes are described as follows. All individuals of type $i$ in either population die at rate 1 and give birth to an individual of type $j$ at rate $\l w_i w_j/n$. 
An offspring of type $j$ at time $t$ is placed in group 1 if and only if the parent was in group 1 and $Z_j^{(1)}(t) = 0$ (that is, there are no individuals of type $j$ in group 1 present already). All other offspring are placed in group 2. 
It is easy to see that $\mathbf{Z}^{(1)}(t)$ has the same distribution as the contact process $\eta_t$. Since $\mathbf{Z}^{(2)}(t) \geq \mathbf{0}$, the proof is complete. 
\end{proof}

Based on the above coupling, we will use $P_\mathbf{w}^{n,i}$ and $P_\mathbf{w}^{n}$ to denote the law of $\mathbf{Z}(t)$ under the initial conditions $\mathbf{Z}(0) = \mathbf{e}_i$ and $\mathbf{Z}(0) = \mathbf{1}$, respectively. 
Part \eqref{dieout} of Theorem \ref{mainthm} follows immediately from Lemma \ref{MTBPcouple} and the following Theorem. 
\begin{thm}
 Let $\l < \l_c$. Then there exists a constant $C>0$ such that
\[
 \lim_{n\ra\infty} P_\mathbf{w}^n (\mathbf{Z}(C \log n) = \mathbf{0} ) = 1, \quad \mu-a.s.
\]
\end{thm}
\begin{proof}
Define the mean matrix $\mathbf{M}(t) = \left( M_{i,j}(t) \right)_{i,j=1}^n$ for the MTBP $\mathbf{Z}(t)$ by
\[
 M_{i,j}(t) = E_\mathbf{w}^{n,i}[ Z_j(t)]. 
\]
Then,
\[
 \frac{d M_{i,j}(t)}{dt} = -M_{i,j}(t) + \sum_{k=1}^n M_{i,k}(t) \frac{\l w_k w_j}{n},
\]
which in matrix form is 
\[
 \frac{d \mathbf{M}(t)}{dt} = -\mathbf{M}(t) \mathbf{A},  
\]
where $\mathbf{A}= (A_{i,j})_{i,j=1}^n$ is the matrix with entries
\[
 A_{i,j} = \frac{\l w_i w_j}{n} - \d_{i,j}.
\]
Solving this system of differential equations gives that 
\[
 \mathbf{M}(t) = e^{\mathbf{A} t}.  
\]
A MTBP such as $\mathbf{Z}(t)$ with mean matrix $\mathbf{M}(t) = e^{\mathbf{A} t}$ is called \emph{sub-critical} if the largest eigenvalue of $\mathbf{A}$ is negative. It is easy to see from the definition of $\mathbf{A}$ above that $\mathbf{A} = \frac{\l}{n}\mathbf{w}_n^* \mathbf{w}_n - I$, where $\mathbf{w}_n = (w_1, w_2, \ldots, w_n)$ is the row vector of weights and $\mathbf{w}_n^*$ is the transpose of $\mathbf{w}_n$. Now, $\mathbf{w}_n^* \mathbf{w}_n$ is a rank one matrix with $n-1$ eigenvalues at $0$ and one eigenvalue at $\| \mathbf{w}_n \|^2 = \sum_{j=1}^n w_j^2$. (To see this, check that $\mathbf{w}_n$ is an eigenvector with this as the corresponding eigenvalue.) 
Therefore, the matrix $\mathbf{A}$ has $n-1$ eigenvalues at $-1$ and one eigenvalue at $-1 + \frac{\l}{n}\| \mathbf{w}_n \|^2$, and the MTBP $\mathbf{Z}(t)$ is sub-critical if 
\[
 \l < \frac{n}{\| \mathbf{w}_n \|^2} = \left( \frac{1}{n} \sum_{j=1}^n w_j^2 \right)^{-1}. 
\]

For any $n$ by $n$ matrix $\mathbf{Q}$, let $\| \mathbf{Q} \| := \sup \{ \| \mathbf{Q} \mathbf{v} \| \,:\, \|\mathbf{v}\| \leq 1 \}$ be the spectral radius of $\mathbf{Q}$. Then, the above analysis of the eigenvalues of $\mathbf{A}$ implies that
\[
 \| \mathbf{M}(t) \| = e^{-t( 1-\frac{\l}{n}\| \mathbf{w}_n \|^2) }. 
\]
The law of large numbers and the definiton of $\l_c$ imply that 
\[
 \lim_{n\ra\infty} \frac{\| \mathbf{w}_n \|^2}{n} = \lim_{n\ra\infty} \frac{1}{n} \sum_{i=1}^n w_i^2 = \frac{1}{\l_c}, \quad \mu-a.s.
\]
Thus, if $\l < \l_c$ there exists an $\e>0$ such that $  1-\frac{\l}{n}\| \mathbf{w}_n \|^2 \geq \e $ for all $n$ large enough, $\mu-a.s$.
Therefore, $\mu-a.s.$,
\[
 P_\mathbf{w}^{n,i}(Z_j(t) \geq 1) \leq E_\mathbf{w}^{n,i}[ Z_j(t) ] = M_{i,j}(t) \leq \| \mathbf{M}(t) \| \leq e^{-\e t}, \quad \text{ for all } n \text{ large enough.}
\]
(The second to last inequality follows from the fact that the entries of $\mathbf{M}(t)$ are non-negative.) 
Therefore, even with the initial configuration of $\mathbf{Z}(0) = \mathbf{1}$ we have that, $\mu-a.s.$,
\[
 P_{\mathbf{w}}^{n} (\mathbf{Z}(t) \neq \mathbf{0} ) \leq \sum_{i,j=1}^n P_\mathbf{w}^{n,i}(Z_j(t) \geq 1) \leq n^2 e^{-\e t}, \quad \text{ for all } n \text{ large enough.} 
\]
Letting $t=C\log n$ for some $C>2/\e$ completes the proof of the Theorem. 
\end{proof}


We can also use the branching process approach to give an upper bound on the probability of a vertex to be infected after a long time. 
The following Proposition complements the lower bound in Corollary \ref{iplb2} and thus completes the proof of Theorem \ref{mainthm2}. 

\begin{prop}\label{rhoub}
 Let $\l > \l_c$. Then, there exists a $C>0$ such that
\[
 \limsup_{n\ra\infty} \sup_{t\geq C\log n} \sup_{i\in[n]} \left( P_\mathbf{w}^n \left( \eta_{t}(i) = 1 \right) - \frac{ \s(\l) w_i}{ 1 + \s(\l) w_i} \right) \leq 0. 
\]
\end{prop}
\begin{proof}
The self-duality property \eqref{duality} and Lemma \ref{MTBPcouple} imply that 
\be\label{etatZt}
 P_{\mathbf{w}}^n ( \eta_t(i) = 1 ) = P_{\mathbf{w}}^{n,i} (\eta_t \neq \mathbf{0} ) \leq P_\mathbf{w}^{n,i} ( \mathbf{Z}(t) \neq \mathbf{0} ).
\ee
To approximate $P_\mathbf{w}^{n,i}( \mathbf{Z}(t) \neq \mathbf{0} ) $ we need the following lemma. 

\begin{lem}\label{extinctionprob}
 For $\l > \left( \frac{1}{n} \sum_{j=1}^n w_j^2 \right)^{-1}$, let 
\[
 \widehat{\s}_n(\l) \text{ be the unique } \s > 0 \text{ that solves } \qquad 1 = \frac{\l}{n} \sum_{j=1}^n \frac{w_j^2}{1+\s w_j}.
\]
Then, the extinction probability $\hat{\rho}_{i,n} = \hat{\rho}_{i,n}(\l)$ satisfies
\[
 \hat{\rho}_{i,n} = \lim_{t\ra\infty} P_\mathbf{w}^{n,i} \left( \mathbf{Z}(t) = \mathbf{0} \right) = \frac{1}{1+\widehat{\s}_n(\l) w_i}. 
\]
\end{lem}
\begin{rem}
If $\l > \l_c$, then $\widehat{\s}_n(\l)$ and $\hat{\rho}_{i,n}$ are defined for all $n$ large enough. Moreover, 
\be\label{rhoidef}
 \lim_{n\ra\infty} \widehat{\s}_n(\l) = \s(\l) \quad\text{and}\quad \lim_{n\ra\infty} \hat{\rho}_{i,n} = \rho_i := \frac{1}{1+\s(\l) w_i}, \qquad \mu-a.s.
\ee
\end{rem}

\begin{proof}
An equation determining the extinction probabilities for multi-type branching processes is given in \cite[Section 7.5]{anBP}. In order to use this, we first need to introduce some notation. Let 
\begin{equation}\label{afdef}
 a_i = 1 + \sum_{j=1}^n \frac{\l w_i w_j}{n}, \qquad f_i(\mathbf{s}) = \frac{1}{a_i} + \sum_{j=1}^n \frac{\l w_i w_j}{a_i n} s_i s_j,  
\end{equation}
and
\begin{equation}\label{udef}
 \mathbf{u}(\mathbf{s}) = ( u_1(\mathbf{s}),u_2(\mathbf{s}),\ldots ,u_n(\mathbf{s})), \quad\text{where}\quad u_i(\mathbf{s}) = a_i( f_i(\mathbf{s}) - s_i ). 
\end{equation}
Then, the extinction probability vector $\hat{\rho}_n = (\hat{\rho}_{1,n}, \hat{\rho}_{2,n}, \ldots,\hat{\rho}_{n,n})$ is the unique solution to 
\begin{equation}\label{uequation}
 \mathbf{u}(\mathbf{s}) = \mathbf{0}, \quad \mathbf{0} \leq \mathbf{s} < \mathbf{1}.
\end{equation}

Let $\s>0$ and $\mathbf{s} = (s_1,s_2,\ldots,s_n)$, where
\[
 s_i = \frac{1}{1 + \s w_i}
\]
Recalling \eqref{afdef} and \eqref{udef} we obtain that
\begin{align*}
 u_i(\mathbf{s}) &= 1 + \frac{\l}{n} w_i s_i \sum_{j=1}^n w_j s_j - s_i - \frac{\l}{n} w_i s_i \sum_{j=1}^n w_j \\
&= (1-s_i) -\frac{\l}{n} w_i s_i \sum_{j=1}^n w_j(1-s_j) \\
&= \frac{\s w_i}{1+\s w_i} - \frac{\l}{n} \frac{w_i}{1+\s w_i} \sum_{j=1}^n \frac{\s w_j^2}{1+\s w_j} \\
&= \frac{\s w_i}{1 + \s w_i}\left( 1 - \frac{\l}{n} \sum_{j=1}^n \frac{ w_j^2}{1 + \s w_j} \right). 
\end{align*}
The proof is completed by noting that the term inside the parenthesis on the last line above equals zero when $\s = \widehat{\s}_n(\l)$. 
\end{proof}
Recall that by \eqref{etatZt} we need an upper bound for $P_\mathbf{w}^{i,n}( \mathbf{Z}(t) \neq \mathbf{0} )$. 
By Lemma \ref{extinctionprob}, we know that for $n$ fixed $P_\mathbf{w}^{i,n}( \mathbf{Z}(t) \neq \mathbf{0} )$ decreases to $ 1-\rho_{i,n}$  as $t$ increases. 
We would like to show that this convergence is fast enough so that the error is very small when $t$ is large enough. 

To this end, let $T_0 := \inf \{ t> 0: \mathbf{Z}(t) = \mathbf{0} \}$ be the extinction time of the MTBP $\mathbf{Z}(t)$. 
Let $\overline{P}_\mathbf{w}^{n,i}(\cdot) = P_\mathbf{w}^{n,i}( \cdot \; | T_0 < \infty)$ be the law of $\mathbf{Z}(t)$ started from one individual of type $i$ and conditioned to eventually die out. Recall that if $\l > \l_c$ then for all $n$ large enough $\hat{\rho}_{i,n} = P_\mathbf{w}^{n,i} ( T_0 < \infty) \in (0,1)$. Then, 
\begin{align}
 P_\mathbf{w}^{n,i}( \mathbf{Z}(t) \neq \mathbf{0} ) &= P_\mathbf{w}^{n,i}( \mathbf{Z}(t) \neq \mathbf{0} \, , \, T_0 < \infty ) + P_\mathbf{w}^{n,i}( \mathbf{Z}(t) \neq \mathbf{0} \, , \, T_0 = \infty ) \nonumber \\
&\leq P_\mathbf{w}^{n,i}( T_0 < \infty ) \overline{P}_\mathbf{w}^{n,i}( \mathbf{Z}(t) \neq \mathbf{0} ) + P_\mathbf{w}^{n,i}( T_0 = \infty )  \nonumber \\
&= \hat{\rho}_{i,n} \overline{P}_\mathbf{w}^{n,i}( \mathbf{Z}(t) \neq \mathbf{0} ) + 1-\hat{\rho}_{i,n}. \label{conditioning}
\end{align}

\begin{lem}\label{cprob}
If $\l > \l_c$ then there exists an $\e>0$ so that, $\mu-a.s.$, for all $n$ large enough
\[
 \overline{P}_\mathbf{w}^{n,i} ( Z_j(t) \neq 0 ) \leq e^{-\e t} , \qquad \forall i,j\leq n. 
\]
\end{lem}
\begin{proof}
The MTBP $\mathbf{Z}(t)$ conditioned to go extinct can be analyzed using $h$-transforms. Given the vertex weights $\mathbf{w}$, let $h_n$ be defined by
\[
 h_n(\mathbf{z}) = \prod_{i=1}^n \hat{\rho}_{i,n}^{z_i}, \qquad \mathbf{z}=(z_1,z_2,\ldots,z_n). 
\]
Then the jump rates of the process $\mathbf{Z}(t)$ conditioned on extinction are given by 
$\bar{q}_n(\mathbf{x},\mathbf{y}) = q_n(\mathbf{x},\mathbf{y})h_n(\mathbf{y})/h_n(\mathbf{x})$, where $q_n(\mathbf{x},\mathbf{y})$ are the jump rates of the process $\mathbf{Z}(t)$ without conditioning on extinction. 
Therefore, we obtain that under $\overline{P}_\mathbf{w}^{n,i}$,  $\mathbf{Z}(t)$ is a MTBP with modified birth and death rates. Under the measuer $\overline{P}_\mathbf{w}^{n,i}$, individuals of type $i$ die at rate $1/\hat{\rho}_{i,n}$ and give birth to particles of type $j$ at rate $\frac{\l w_i w_j}{n} \hat{\rho}_{j,n}$. 
We define the mean matrix conditioned on extinction by $\overline{\mathbf{M}}(t) = (\overline{M}_{i,j}(t) )_{i,j=1}^n$, where 
\[
 \overline{M}_{i,j}(t) = \overline{E}_{\mathbf{w}}^{n,i} [ Z_j(t)],
\]
and $\overline{E}_{\mathbf{w}}^{n,i}$ denotes expectations with respect to $\overline{P}_\mathbf{w}^{n,i}$. 
Then, as was done for the mean matrix $\mathbf{M}(t)$ above we see that 
\[
 \overline{\mathbf{M}}(t) = e^{\bar{\mathbf{A}} t}, \quad\text{where}\quad \bar{\mathbf{A}} = (\bar{A}_{i,j})_{i,j=1}^n, \quad 
\bar{A}_{i,j} = 
\begin{cases}
 -\frac{1}{\hat{\rho}_{j,n}} + \frac{\l w_i w_j}{n} \hat{\rho}_{j,n} & i=j  \\
\frac{\l w_i w_j}{n} \hat{\rho}_{j,n} & i\neq j.
\end{cases}
\]
Letting $D_{\hat{\rho}_n}$ be the diagonal matrix with diagonal entries $\hat{\rho}_{i,n}$ we obtain that
\[
\overline{ \mathbf{A}} = \frac{\l}{n} \mathbf{w}_n \mathbf{w}_n^* D_{\hat{\rho}_n} - D_{\hat{\rho}_n}^{-1} = D_{\hat{\rho}_n}^{-1} \left( \frac{\l}{n} D_{\hat{\rho}_n} \mathbf{w}_n \mathbf{w}_n^* D_{\hat{\rho}_n} - I \right)
\]
As above, we see that the eigen values of $\frac{\l}{n} D_{\hat{\rho}_n} \mathbf{w}_n \mathbf{w}_n^* D_{\hat{\rho}_n} - I $ are $-1$ (with multiplicity $n-1$) and $\frac{\l}{n} \| D_{\hat{\rho}_n} \mathbf{w}_n \|^2 - 1 = \frac{1}{n} \sum_{i=1}^n w_i^2 \hat{\rho}_{i,n}^2 - 1$. Note that the defnitions of $\hat{\rho}_{i,n}$ and $\widehat{\s}_n(\l)$ imply that 
\[
 \frac{1}{n} \sum_{i=1}^n w_i^2 \hat{\rho}_{i,n} = \frac{1}{n} \sum_{i=1}^n \frac{w_i^2}{1+\widehat{\s}_n(\l)}  = 1. 
\]
Therefore, since $\hat{\rho}_{i,n} < 1$ we obtain that (for $n$ large enough) $\frac{1}{n} \sum_{i=1}^n w_i^2 \hat{\rho}_{i,n}^2 < 1$, and thus all the eigenvalues of $\frac{\l}{n} D_{\hat{\rho}_n} \mathbf{w}_n \mathbf{w}_n^* D_{\hat{\rho}_n} - I $ are negative. We wish, however, to show that all the eigenvalues of $\overline{\mathbf{A}}$ are negative. To this end, note that $D_{\hat{\rho}_n}^{-1}$ is a diagonal matrix with diagonal entries all greater than one. We prove the following algebra lemma in the appendix. 
\begin{lem}\label{eigenvalues}
 Let $U$ be a symmetric matrix with all eigenvalues less than $\l_0 < 0$, and let $D$ be a diagonal matrix with diagonal entries all greater than $1$. Then, $DU$ also has all eigenvalues less than $\l_0$. 
\end{lem}
Applying Lemma \ref{eigenvalues}, we obtain that $\overline{\mathbf{A}}$ has all eigenvalues less than $\frac{1}{n} \sum_{i=1}^n w_i^2 \hat{\rho}_{i,n}^2 - 1 < 0$. 
Recalling the definition of $\rho_i$ in \eqref{rhoidef}, we obtain that
\[
 \lim_{n\ra\infty} \frac{1}{n} \sum_{i=1}^n w_i^2 \hat{\rho}_{i,n}^2 - 1 = E_\mu\left[ w_1^2 \rho_1^2 \right] - 1 = E_\mu\left[\frac{w_1^2}{(1+\s(\l) w_1)^2} \right] - 1 < 0, 
\]
and therefore, $\mu-a.s.$, for any $0 < \e < 1 - E_\mu\left[ w_1^2/(1+\s(\l) w_1)^2 \right]$ and all $n$ large enough, $\| \overline{\mathbf{M}}(t) \| \leq e^{-\e t}$. 
The proof is completed by noting that 
\[
 \overline{P}_\mathbf{w}^{n,i} ( Z_j(t) \neq 0 ) \leq \overline{M}_{i,j}(t) \leq \| \overline{\mathbf{M}}(t) \|. 
\]
\end{proof}
We are now ready to finish the proof of Proposition \ref{rhoub}. Recalling \eqref{etatZt} and \eqref{conditioning}, Lemmas \ref{extinctionprob} and \ref{cprob} imply that, $\mu-a.s.$, for all $n$ large enough,
\[
P_\mathbf{w}^n (\eta_t(i) = 1)
\leq n e^{-\e t} + 1-\hat{\rho}_{i,n} = n e^{-\e t} + \frac{\widehat{\s}_n(\l) w_i}{1+ \widehat{\s}_n(\l) w_i} , \qquad \forall i\leq n. 
\]
Therefore, recalling \eqref{etatZt} we obtain that
\begin{align*}
& \sup_{t\geq C\log n} \sup_{i\in[n]} \left( P_\mathbf{w}^n (\eta_t(i) = 1) - \frac{\s(\l) w_i}{1+\s(\l) w_i} \right) \\
&\qquad \leq \sup_{t\geq C\log n} \sup_{i\in[n]} \left( n e^{-\e t} + \frac{\widehat{\s}_n(\l) w_i}{1+ \widehat{\s}_n(\l) w_i} - \frac{\s(\l) w_i}{1+\s(\l) w_i} \right) \\
&\qquad \leq n^{1-C\e}  + \frac{| \widehat{\s}_n(\l) - \s(\l)|}{\left( \sqrt{\widehat{\s}_n(\l)}  + \sqrt{\s(\l)}  \right)^2},
\end{align*}
where we used \eqref{calculus} in the last inequality. 
Since $\widehat{\s}_n(\l) \ra \s(\l)$ as $n\ra\infty$, the conclusion of Proposition \ref{rhoub} holds for any $C>1/\e$. 
\end{proof}

\section{Asymptotics of $\s(\l)$}\label{criticalexponent}

In this section we give the proof of Proposition \ref{cexp}. 
The proof doesn't involve probability at all and is purely an exercise in analysis.

Recall that $\s(\l)$ is defined for $\l > \l_c$ by 
\[
 \s(\l) \text{ is the unique } \s > 0 \text{ that solves } \qquad 1 = \l E_\mu\left[ \frac{ w_1^2}{1+ \s w_1} \right]. 
\]
Thus, it is crucial to understand the asymptotics of $E_\mu\left[  \frac{ w_1^2}{1+ \s w_1} \right]$ as $\s\ra 0^+$. 
\begin{lem}\label{Easymp}
 If $\mu(w_1 \geq x) \sim C x^{-(\a-1)}$ for some $C>0$ and $\a>2$, then 
\be \label{a23}
 E_\mu\left[  \frac{ w_1^2}{1+ \s w_1} \right] \sim
\begin{cases}
 C \frac{(\a-1)\pi}{\sin(\pi \a)} \s^{\a-3} & \a \in (2,3) \\
 -2C \log \s & \a = 3, 
\end{cases}
\qquad \text{ as } \s \ra 0^+,
\ee
and 
\be\label{a34}
 E_\mu\left[  \frac{ w_1^2}{1+ \s w_1} \right] - \frac{1}{\l_c} \sim 
\begin{cases}
 C \frac{(\a-1)\pi}{\sin(\pi \a)} \s^{\a-3} & \a \in (3,4) \\
 3 C \s \log \s & \a = 4 \\
 - E_\mu[w_1^3] \s & \a > 4, 
\end{cases}
\qquad \text{ as } \s \ra 0^+.
\ee
Moreover, the conclusion in the case $\a>4$ also holds under the assumption that $E_\mu[w_1^3] < \infty$. 
\end{lem}
\begin{proof}
 We first consider the case when $\a \leq 3$. Note that in this case $E_\mu[w_1^2] = \infty$, and thus $E_\mu[w_1^2/(1+\s w_1)] \ra \infty$ as $\s\ra 0^+$. The assumption on the tail decay of $\mu$ implies that for any $\e>0$, there exists an $x_0= x_0(\e) < \infty$ such that 
\[
 (C-\e) x^{-(\a-1)} \leq \mu( w_1 \geq x) \leq (C+\e) x^{-(\a-1)}, \qquad \forall x \geq x_0.
\]
Therefore,
\begin{align*}
 \int_{x_0}^{\infty} \frac{x^2}{1+\s x} \left( (C-\e)(\a-1) x^{-\a} \right) dx & \leq E_\mu\left[  \frac{ w_1^2}{1+ \s w_1} \right]\\
& \leq \frac{x_0^2}{1+\s x_0} + \int_{x_0}^{\infty} \frac{x^2}{1+\s x} \left( (C+\e)(\a-1) x^{-\a} \right) dx.
\end{align*}
Note that $\frac{x_0^2}{1+\s x_0}$ and $\int_1^{x_0}\frac{x^{2-\a}}{1+\s x} dx$ both are bounded as $\s\ra 0^+$. Therefore it is enough to show
\be\label{intasymp}
 \int_{1}^{\infty} \frac{x^{2-\a}}{1+\s x} dx \sim
\begin{cases}
 \frac{\pi}{\sin(\pi \a)} \s^{\a-3} & \a \in (2,3) \\
 -\log \s & \a = 3, 
\end{cases}
\qquad \text{ as } \s \ra 0^+.
\ee
When $\a=3$, we can explicitly evaluate this integral to be $\log(1+1/\s)$, which is asymptotic to $ -\log \s$ as $\s\ra 0^+$. 
To handle the case when $\a\in (2,3)$ we first make the change of variables $z= \s x$ in the integral in \eqref{intasymp} to obtain
\[
 \int_{1}^{\infty} \frac{x^{2-\a}}{1+\s x} dx = \s^{\a-3} \int_{\s}^\infty \frac{z^{2-\a}}{1+z} dz. 
\]
However, taking $\s\ra 0^+$ in the limit of the last integral we obtain that for any $\a\in(2,3)$
\[
 \lim_{\s \ra 0^+} \int_{\s}^\infty \frac{z^{2-\a}}{1+z} dz = \int_{0}^\infty \frac{z^{2-\a}}{1+z} dz = \frac{\pi}{\sin(\pi \a)}, 
\]
where the last equality is a standard exercise in complex analysis.
This completes the proof of \eqref{intasymp} and thus also the proof of \eqref{a23}. 

Next, we consider the case when $\a>3$. In this case $E_\mu[w_1^2] = 1/\l_c < \infty$, and thus
\be\label{soasymp}
 E_\mu\left[  \frac{ w_1^2}{1+ \s w_1} \right] = E_\mu\left[  w_1^2\left(1-\frac{\s w_1}{1+\s w_1} \right) \right] 
= \frac{1}{\l_c} - \s E_\mu \left[\frac{w_1^3}{1+\s w_1} \right].
\ee
Thus, we are reduced to computing the asymptotics of $E_\mu[w_1^3/(1+\s w_1)]$. 
Note that $E_\mu[w_1^3/(1+\s w_1)] \ra E_\mu[w_1^3]$ as $\s\ra 0^+$. 
This completes the proof of \eqref{a34} in the case when $E_\mu[w_1^3] < \infty$ (in particular when $\a>4$). 
On the other hand, if $\a \leq 4$ then $E_\mu[w_1^3] = \infty$ and we must compute the asymptotics of how fast $E_\mu[w_1^3/(1+\s w_1)]$ grows as $\s\ra 0^+$. 
To this end, we repeat the process above in the case when $\a\in (2,3]$ to show that 
\[
 E_\mu \left[\frac{w_1^3}{1+\s w_1} \right] 
\sim
 C (\a-1) \int_1^\infty \frac{x^{3-\a}}{1+\s x} dx
\sim 
\begin{cases}
 - C \frac{(\a-1)\pi}{\sin(\pi \a)} \s^{\a-4} & \a \in (3,4) \\
 -3 C \log \s & \a = 4,
\end{cases}
\quad \text{ as } \s \ra 0^+.
\]
Combining these asymptotics with \eqref{soasymp} completes the proof of \eqref{a34} when $\a \in (3,4]$. 
\end{proof}

With Lemma \ref{Easymp} at our disposal we are now ready to give the proof of Proposition \ref{cexp}. 
\begin{proof}[Proof of Proposition \ref{cexp}]
 First, note that the definition of $\s(\l)$ and the asymptotics in Lemma \ref{Easymp} imply that $\s(\l) \ra 0$ as $\l \ra\l_c$. Then, another application of Lemma \ref{Easymp} implies that
\[
 1 = \l E_\mu\left[ \frac{ w_1^2}{1+ \s(\l) w_1} \right] 
\sim 
\begin{cases}
 \l C \frac{(\a-1)\pi}{\sin(\pi \a)} \s(\l)^{\a-3} & \a \in (2,3) \\
 -2C \l \log \s(\l) & \a = 3, 
\end{cases}
\qquad \text{ as } \l \ra 0^+. 
\]
The above asymptotics are equivalent to the conclusion of Proposition \ref{cexp} in the cases $\a\in(2,3)$ and $\a=3$. 
Similarly, when $\a>3$ the defenition of $\s(\l)$ and Lemma \ref{Easymp} imply that 
\[
\frac{\l-\l_c}{\l_c} = \l\left( \frac{1}{\l_c} - E_\mu \left[ \frac{w_1^2}{1+\s(\l) w_1} \right] \right)
\sim
\begin{cases}
  - \l C \frac{(\a-1)\pi}{\sin(\pi \a)} \s(\l)^{\a-3} & \a \in (3,4) \\
 - \l 3 C \s(\l) \log \s(\l) & \a = 4 \\
 \l E_\mu[w_1^3] \s(\l) & \a > 4, 
\end{cases}
\qquad \text{ as } \l \ra \l_c^+. 
\]
The conclusions of Proposition \ref{cexp} follow easily from the above asymptotics in the cases $\a\in(3,4)$ and $a>4$ (or $E_\mu[ w_1^3]<\infty$). 
In the case $\a=4$, the conclusion in Proposition \ref{cexp} follows from the above asymptotics and the fact that 
\[
 f(\d)\log f(\d) \sim -A \d  \iff f(\d) \sim A \frac{\d}{\log(1/\d)},
\]
where the above asymptotics are as $\d \ra 0^+$. 
\end{proof}

\appendix

\section{Proof of Lemma \ref{eigenvalues}}
 
The proof of Lemma \ref{eigenvalues} is a rather simple linear exercise in linear algebra, but we give the proof here for completeness.
\begin{proof}
As in the statement of the Lemma, let $U$ be a real $n\times n$ symmetric matrix with largest eigenvalue $\l_0 < 0$, and let $D$ be an $n\times n$ diagonal matrix with diagonal entries $D_{i,i} \geq 1$. First, note that $DU$ has the same eigenvalues as $D^{-1/2}DUD^{1/2} = D^{1/2}UD^{1/2}$. Note that since the latter is obviously symmetric, $DU$ has all real eigenvalues. 

Let $\theta_0$ be the largest eigenvalue of $DU$ (and also of $D^{1/2}UD^{1/2}$). 
Then, by the well known variational characterisation of the largest eigenvalue of symmetric matrices
\begin{align*}
 \theta_0 &= \sup \{ \langle D^{1/2} U D^{1/2} u, u \rangle : \, |u|=1 \} \\
&= \sup \{  \langle  U D^{1/2} u, D^{1/2} u \rangle : \, |u|=1 \} \\
&\leq \sup \{  \langle  U v, v \rangle : \, |v|\geq 1 \} \\
&= \sup \{  \langle  U v, v \rangle : \, |v| = 1 \} = \l_0. 
\end{align*}
Note that the inequality above follows from the fact that $| D^{1/2} u | \geq |u|$, and the second to last equality holds since $U$ is negative-definite. 
\end{proof}

\bibliographystyle{alpha}
\bibliography{RandomGraphs}

\end{document}